\documentclass[11pt,a4paper]{article}
\usepackage{amsmath,amssymb}

\usepackage{array,longtable}
\usepackage{color}
\usepackage[colorlinks,linkcolor=blue,citecolor=red]{hyperref}
\usepackage{graphicx}
\usepackage{wrapfig}

\DeclareMathSymbol{\bbbr}{\mathalpha}{AMSb}{"52}
\DeclareMathSymbol{\bbbc}{\mathalpha}{AMSb}{"52}

\newtheorem{theorem}{Theorem}

\newtheorem{corollary}[theorem]{Corollary}

\newtheorem{definition}{Definition}

\newtheorem{lemma}[theorem]{Lemma}

\newtheorem{proposition}[theorem]{Proposition}

\begin{document}

\title{Fractional-linear integrals of geodesic flows on surfaces and  Nakai's geodesic 4-webs}

\author{{\Large Sergey I. Agafonov}\\
\\
{\Large Thaís G. P. Alves}\\
\\
Department of Mathematics,\\
S\~ao Paulo State University-UNESP,\\ S\~ao Jos\'e do Rio Preto, Brazil\\
e-mail: {\tt sergey.agafonov@gmail.com} \\
 \ \ \ \ \ {\tt thaguinami@hotmail.com} }
\date{}
\maketitle
\unitlength=1mm

\vspace{1cm}

\begin{abstract}
We prove that if the geodesic flow on a surface has an integral,  fractional-linear in momenta, then the dimension of the space of such integrals is either  3 or 5, the latter case corresponding to constant gaussian curvature. We give also a geometric criterion for existence of fractional-linear integrals: such integral exists if and only if the surface carries a geodesic 4-web with constant cross-ratio of the four directions tangent to the web leaves. 
\end{abstract}

\section{Introduction} 
Integrable geodesic flow on surfaces is a classical subject of differential geometry and analytic mechanics since the 19th century works of Jacobi \cite{J-84}, Darboux \cite{D-91}, Dini \cite{D-69}, Koenigs \cite{K-96} and many others. Special attention was payed to  integrals  polynomial in momenta. This problem is well motivated even in the local setting due to the following observations: 1) if there is a locally defined  analytic integral then each  homogeneous in momenta part of its Taylor expansion is an integral  polynomial in momenta, 2) the surface of zero velocities  in the tangent bundle of surface is the set of singular points of the geodesic flow, the existence of regular integrals at singular point being a nontrivial restriction. 

Since the geodesic flow is Hamiltonian,   it is enough to find one first integral independent of the Hamiltonian function to integrate the flow due to the Liouville Theorem. Polynomial integrals of fixed degree constitute a vector space of finite dimension. Establishing this dimension turned out to be a non-trivial task for degrees grater then two, in fact, the list of possible dimensions is not known already for cubic integrals (see \cite{K-08}). The case of linear integrals is equivalent to listing the possible dimensions of isometry group due to the Noether Theorem: the dimension is either 0, or 1, or 3. The case of quadratic integrals was settled by Koenings: the possible dimensions are 1, 2, 3, 4, and 6, the largest dimension in both cases corresponding to constant gaussian curvature. 

Modern interests to this topic is due to 1) a new discovered relation to infinite-dimensional integrable systems of PDE and 2) a better understanding of geometry behind the integrability of geodesic flows on surfaces.  

For a fixed metric, the system of partial differential equations (PDE) for the coefficients of the integral  is over-determined but if one considers the metric also as unknown then the system turns out to be of hydrodynamic type with remarkable properties: it is diagonalizable,  it possesses infinitely many conservation laws, and it is linearizable by Tsarev's {\it generalized hodograph method} (see \cite{K-82,T-97,BM-11,MP-17,T-85}). 

On the other hand, the existence of quadratic and cubic integrals can be described in geometric terms via geodesic webs. For example, consider a cubic integral, which can be rewritten also as polynomial in velocity $v=(\xi,\eta)$ due to the canonical isomorphism between the tangent and cotangent bundles of the surface with local coordinates $x,y$:  
$$
I_3=a_0(x,y)\xi^3+a_1(x,y)\xi^2\eta+a_2(x,y)\xi\eta^2+a_3(x,y)\eta^3.
$$
Equation it to zero, one gets an implicit cubic ordinary differential equation (ODE). If all its 3 roots $[\xi:\eta]$  are real, the integral curves of this ODE form a hexagonal geodesic 3-web. And conversely, the existence of  a hexagonal geodesic 3-web implies existence of a cubic integral \cite{Ah-21}.   
There is also a geometric characterization of quadratic integrals \cite{Aq-21}:   existence of quadratic integrals is equivalent to existence of geodesic net  $\mathcal{G}$ such that any  3-subweb of 4-web, formed by  $\mathcal{G}$ and its bisector net  $\mathcal{N}$, is hexagonal. 

Geodesic webs is  a classical chapter of differential geometry, such webs were a subject of intensive study, though without any relation to dynamics, see \cite{F-99,S-03,GS-24,S-26,V-29}.  In his last book on the web theory \cite{B-55} Blaschke posed the problem of finding an intrinsic criterion for existence of hexagonal geodesic 3-webs. Thus, the result of \cite{Ah-21} provides such criterion it dynamical terms and generalizes the Graf and Sauer Theorem \cite{GS-24} for metrics of non-constant gaussian curvature.      

Another natural class of integrals of geodesic flows are integrals rational in momenta, whose numerator and denominator are homogeneous polynomials in momenta of the same degree. The study of such integrals was initiated by Kozlov \cite{K-14}. 

Agapov and Shubin investigated fractional-linear integrals \cite{Ar-20,At-20,AS-21}    
$$\dfrac{A(x,y)p+B(x,y)q}{C(x,y)p+E(x,y)q},$$
in particular, they found explicit examples of metrics, whose geodesic flows admit such integrals.

This paper contributes to both  analytic and geometric description of fractional-linear integrals. The Möbius group acts naturally on the space of such integrals,  thus giving a symmetry subgroup of the PDE system for the integral. Splitting this system to {\it automorphic} and {\it resolving} parts (see \cite{V-04,O-82}), we manage to compute enough compatibility conditions to extract the dimension of the space of fractional-linear integrals.  

The first main result of the paper claims that if this space is nonempty then  this dimension is either 3 or 5, the maximal dimension corresponding to constant gaussian curvature. 

To give a geometric description, we construct a one-parametric family of geodesic foliations, whose leaves are integral curves of the ODE 
$$\dfrac{\alpha(x,y)\xi+\beta(x,y)\eta}{\gamma(x,y)\xi+\delta(x,y)\eta}=\lambda,$$
where on the left-hand side is a fractional-linear integral written on the tangent bundle and $\lambda \in \mathbb{R}$ is a parameter. Taking four such foliations one obtains a geodesic 4-web with constant cross-ratio of web tangent directions.  

Nakai studied 4-webs with constant cross-ratio of web tangent directions in \cite{N-98} and showed that the Blaschke curvatures of all its 3-subwebs are equal. He excluded the case of parallelizable 4-web with hexagonal 3-subwebs as trivial, though we will tolerate this degeneracy in our definition, as it occurs for metrics with constant gaussian curvatures. 
 Dufour e Lehmann proved  \cite{DL-22} that Nakai's webs can have rank  0 or 1 and gave examples of webs with rank 1.
 
Our second main results claims that there is a fractional linear integral if and only if the surface carries Nakai's geodesic 4-web. 

All the results of this paper are local, all functions, surfaces, fields, and other objects are smooth. 

\section{Fractional-linear integrals of geodesic flow}
Geodesic flow on a surface  $S$ with local coordinates $(x,y)$ is Hamiltonian   with respect to canonical symplectic form $\Omega= dp \wedge dx+dq \wedge dy$ on the cotangent bundle $T^*S$, where $f=(p,q)$ is the momentum. The Hamiltonian  $\frac{1}{2}g(v,v)$, where $v=(\xi,\eta)\in T_{s}S$ is the velocity and $g$ metric, can be rewritten in momentum due to the canonical isomorphism between $T_{s}S$ e $T_{s}^*S$ via $g(v,*)=\langle f,*\rangle$.    
For metric $g=\Lambda(x,y)(dx^2+dy^2)$ in conformal coordinates, we get the Hamiltonian $H=\dfrac{p^2+q^2}{2\Lambda}$ and $p,q$ are conjugated to the coordinates $(x,y)$.  

\begin{lemma}\cite{K-14}\label{lema sistema homo}
Geodesic flow of the metric  $g=\Lambda(x,y)(dx^2+dy^2)$ admits a fractional-linear integral \begin{equation}\label{fl}
	I=\frac{A(x,y)p+B(x,y)q}{C(x,y)p+E(x,y)q}
\end{equation} if and only if 
\begin{equation}\label{inteqconform}
 \begin{array}{l}
			2\Lambda( CA_x-AC_x) - \Lambda_y(AE - BC)=0,\\
			2\Lambda(EA_x-AE_x+CA_y-AC_y+CB_x-BC_x) + \Lambda_x(AE - BC)=0,\\
			2 \Lambda (EA_y-AE_y+CB_y-BC_y+EB_x-BE_x) - \Lambda_y(AE - BC)=0,\\
			2\Lambda(EB_y-BE_y) + \Lambda_x(AE - BC)=0.
		\end{array}
\end{equation}
	
\end{lemma}
{\it Proof:} 
The condition $\{I,H\}=0$ writes as an equation, cubic and  homogeneous  in $p,q$:
 \begin{eqnarray*}
 	& &(2\Lambda CA_x-2\Lambda AC_x  -\Lambda_y AE  + \Lambda_y BC)p^3\\
 	 &+& (2\Lambda EA_x- 2\Lambda AE_x+2\Lambda CA_y-2\Lambda AC_y+ 2 \Lambda CB_x - 2\Lambda BC_x +\Lambda_x AE 	 
 	-\Lambda_x BC )p^2q\\ 
 	&+& ( 2 \Lambda EA_y- 2 \Lambda AE_y+ 2 \Lambda CB_y- 2 \Lambda BC_y+ 2 \Lambda EB_x- 2 \Lambda BE_x-\Lambda_y AE  
 	 -\Lambda_y BC  )pq^2\\ 
 	 &+& (2 \Lambda EB_y- 2 \Lambda BE_y+\Lambda_x AE - \Lambda_x BC)q^3=0.
 \end{eqnarray*}
Since $A, B, C$ e $E$ are independent of  $p, q$, the equation splits into four given in (\ref{inteqconform}).
\hfill $\Box$\\

  \section{Differential invariants}
Let us normalize a fractional linear integral (\ref{fl}) to $A(x,y)E(x,y)-B(x,y)C(x,y)=1.$\\ 
Thus normalized  integral defines the map
\begin{equation}\label{mmap}
 \begin{array}{l}
 m: S \longrightarrow SL_2(\mathbb{R}), \\
 \\
 (x,y) \mapsto m(x,y)=\left(\begin{array}{cc} A(x,y) & B(x,y)  \\ C(x,y) & E(x,y)\end{array}\right).
 \end{array}
\end{equation}
The matrix Lie group $SL_2(\mathbb{R})=\left\{\left(\begin{array}{c c}
	\alpha & \beta\\
	\gamma & \delta
\end{array}\right): \ \alpha, \beta, \gamma, \delta \in \mathbb{R} , \alpha\delta-\beta\gamma=1 \right\}$ naturally acts on the set  of fractional-linear integrals:

\begin{equation}
h \cdot I= \frac{\alpha I+\beta}{\gamma I+\delta}=\dfrac{(\alpha A+\beta C)p+(\alpha B+\beta E)q}{(\gamma A+\delta C)p+(\gamma B+\delta E)q},
\end{equation}
where $h= \left(\begin{array}{c c}
	\alpha & \beta\\
	\gamma & \delta
\end{array}\right)\in SL_2(\mathbb{R})$. Note that 
$$\tilde{I}=\dfrac{(\alpha A+\beta C)p+(\alpha B+\beta E)q}{(\gamma A+\delta C)p+(\gamma B+\delta E)q}$$ 
is also normalized 
$$(\alpha A+\beta C)(\gamma B+\delta E)-(\alpha B+\beta E)(\gamma A+\delta C)=1.$$
\begin{definition}
	We say that two fractional-linear integrals  $I$ e $\tilde{I}$ are  equivalent  if there exists an element $h \in SL_2(\mathbb{R})$ such  that  $h\cdot I = \tilde{I}$. 
\end{definition}
Thus the system (\ref{inteqconform}) is invariant by the above defined action of $SL_2(\mathbb R)$ and can be rewritten in differential invariants.  Since the action is just multiplying $m$ by $h$ on the left,  differential invariants can be obtained via {\it Darboux derivative} of the map $m$ defined by (\ref{mmap}). Recall that the Darboux derivative of the map 
 $m:S \rightarrow SL_2(\mathbb{R})$ is the pull-back of the Maurer-Cartan form of $SL_2(\mathbb{R})$: $\omega=m^{-1}\cdot dm$ (see \cite{S-96}).
 In fact, it is invariant as 
 $$(hm)^{-1}\cdot d(hm)=m^{-1}\cdot h^{-1}\cdot h\cdot d(m)=m^{-1}\cdot d(m).$$  We have
\begin{equation}\label{eq3}
	m^{-1}\cdot dm=\left(\begin{array}{cc} E & -B  \\ -C & A\end{array}\right)\left(\begin{array}{cc} dA & dB  \\ dC & dE\end{array}\right)=\left(\begin{array}{cc} EdA-BdC & EdB-BdE  \\ -CdA+AdC & -CdB+AdE\end{array}\right), 
\end{equation}
where 
$$dA=A_xdx+A_ydy,\ \ \ dB=B_xdx+B_ydy,\ \ \ dC=C_xdx+C_ydy,\ \ \  dE=E_xdx+E_ydy.$$
Substituting into (\ref{eq3}) we obtain

	\begin{eqnarray*} \label{eq4}
m^{-1}\cdot dm=
	\left(\begin{array}{cc} EA_x-BC_x & EB_x-BE_x  \\ AC_x -CA_x&AE_x -CB_x\end{array}\right)dx+\left(\begin{array}{cc} EA_y-BC_y & EB_y-BE_y  \\ AC_y-CA_y& AE_y-CB_y \end{array}\right)dy.
\end{eqnarray*}
The group does not act on independent coordinates $x,y$, therefore all the elements of the matrix coefficients of $dx$ and $dy$ are scalar diferential invariants. Since the trace of $\omega$ vanishes we have 6 invariants
$$ 
\begin{array}{lll}
P:=EA_x-BC_x,& Q:=EB_x-BE_x,& R:=AC_x-CA_x,\\
&&\\
X:=EA_y-BC_y, & Y:=EB_y-BE_y, & Z:=AC_y-CA_y.
\end{array}
$$
Using the normalization  $det(m)=1$ we exclude $E=\frac{BC+1}{A}$ and its derivatives
\begin{equation}\label{eq5}
	E_x=\dfrac{BC_x+CB_x}{A} - \dfrac{(BC + 1)A_x}{A^2},\ \ \ \ 
	E_y= \dfrac{BC_y+CB_y}{A} - \dfrac{(BC + 1)A_y }{A^2}, 
\end{equation} 
and rewrite the invariants
$$ 
\begin{array}{lll}
P=\frac{BC + 1}{A}A_x - BC_x,& Q=\dfrac{B^2C+B}{A^2}A_x+\dfrac{1}{A}B_x-\dfrac{B^2}{A}C_x,& R=AC_x-CA_x,\\
&&\\
X=\frac{BC + 1}{A}A_y - BC_y, & Y=\dfrac{B^2C+B}{A^2}A_y+\dfrac{1}{A}B_y-\dfrac{B^2}{A}C_y, & Z=AC_y-CA_y.
\end{array}
$$
We can express   derivatives  $A_x,B_x,C_x,A_y,B_y,C_y$ in terms of the invariants as follows:
\begin{equation} \label{derivadas parciais}
\begin{array}{lll}
A_x=PA+RB, & B_x=QA-PB, & C_x=\dfrac{PAC + RBC + R}{A},\\
&&\\
A_y=XA+ZB, & B_y=YA - XB, & C_y=\dfrac{XAC + ZBC + Z}{A}.\\
\end{array}
\end{equation}
Calculation of compatibility conditions is greatly simplified for the metric form 
\begin{equation}\label{complex}
g=L(x,y)dxdy,
\end{equation}
 related by formal complex substitution $x=u+iv$, $y=u-iv$ to the conformal form:
$$g=L(x,y)dxdy=\Lambda(u+iv,u-iv)(du^2+dv^2).$$
Since all the expressions are rational, this substitution does not affect our compatibility results obtained in a purely algebraic way  (compare with calculations of Lie in  \cite{L-82} and Koenigs in \cite{K-96}).

\begin{lemma}
Geodesic flow on a surface admits a fractional-linear integral (\ref{fl}) if and only if  the invariants $P,Q,R,X,Y,Z$ are related to the metric  (\ref{complex}) as follows:
	\begin{equation}\label{inv}
P=-\dfrac{Y}{2}-\dfrac{L_x}{2L},\,\ \ X=\dfrac{R}{2}+\dfrac{L_y}{2L},\ \ Q=Z=0.
\end{equation}
\end{lemma}
{\it Proof:} 
For the metric (\ref{complex}) the  Hamiltonian assumes the form $H=\dfrac{2pq}{L},$ then the condition  $\{I,H\}=0,$ writes as 
$$
	\begin{array}{l}
		L(AC_y - CA_y)=0,\\
		L(CA_x-AE_y-AC_x+CB_y-BC_y+EA_y) - L_y(AE - BC)=0,\\
		L(EA_x-BE_y-BC_x+EB_y-AE_x+CB_x)+ L_x(AE - BC)=0,\\
		L(BE_x - EB_x)=0.
	\end{array}
$$                                                                   
Using the normalization  $AE-BC=1$, substituting $E,E_x,E_y$ from (\ref{eq5}) and  $A_x,B_x,C_x,A_y,B_y,C_y$ from (\ref{derivadas parciais}) we get 
$$
	LZ=LQ=0,\ \ 2LX-LR-L_y=0,\ \ 2LP+LY+L_x=0,\\
$$
hence (\ref{inv}).  
\hfill $\Box$\\

\section{Dimension of the  space of integrals}
The gaussian curvature of metric (\ref{complex}) is 
	\begin{equation}\label{curvatura gaussiana}
	K=-\dfrac{LL_{xy}-L_xL_y}{L^3},
	\end{equation} 
its partial derivatives are 
 \begin{eqnarray*}
 	K_x=-\dfrac{L_{xxy}}{L^2}+\dfrac{L_{xx}L_y+3L_{xy}L_x}{L^3}-\dfrac{3L_x^2L_y}{L^4},\\ \\
 		K_y=-\dfrac{L_{xyy}}{L^2}+\dfrac{L_{yy}L_x+3L_{xy}L_y}{L^3}-\dfrac{3L_y^2L_x}{L^4}.
 \end{eqnarray*}
The Darboux derivative $\omega$ of $m$ satisfies the structure equation $d\omega+\omega\wedge\omega=0$, which is the local condition for existence of $m$ with the Darboux derivative $\omega$ (see \cite{S-96}). The matrix structure equation rewrites as three scalar equations
$$
\begin{array}{l}
	2RY+2LK - R_x - Y_y =0,\\ 
	 \\
	LY^2+L_xY-LY_x=0,\\ 
	\\
	LR^2+L_yR-LR_y=0.\\ 
\end{array}
$$
Resolving for  $Y_x$ e $R_y$ we obtain 
$Y_x= Y^2+\dfrac{L_x}{L}$ and 
$R_y =R^2+ \dfrac{L_y}{L}.$
Let us introduce $F$ by
\begin{equation}\label{eq15}
	2F:=R_x-Y_y,
\end{equation}
then 
$R_x=RY+LK+F, \ \ Y_y=RY+LK-F.$ 
Thinking of $F$ as an unknown function of $x,y$, we 
get the following exterior differential system in  $\mathbb{R}^8$ with coordinates $A, B, C, R, Y, F, x, y$: 
\begin{eqnarray}\ \label{dA}
	dA&=&\left(RB-\dfrac{L_x}{2L}A-\dfrac{YA}{2}\right)dx+\left(\dfrac{L_y}{2L}A+\dfrac{RA}{2}\right)dy,\\ \nonumber \\ \label{dB}
	dB&=&\left(\dfrac{L_x}{2L}B+\dfrac{YB}{2}\right)dx+\left(YA-\dfrac{L_y}{2L}B-\dfrac{RB}{2}\right)dy,\\\nonumber \\ \label{dC}
	dC&=&\left(\dfrac{BRC}{A}+\dfrac{R}{A}-\dfrac{L_x}{2L}C-\dfrac{YC}{2}\right) dx+ \left(\dfrac{L_y}{2L}C+\dfrac{RC}{2}\right)dy,\\\nonumber \\ \label{dR}
	dR&=&(YR+LK+F)dx+\left(R^2+\dfrac{L_y}{L}R\right)dy,\\\nonumber \\ \label{dY}
	dY&=&\left(Y^2+\dfrac{L_x}{L}Y\right)dx+(RY+LK-F)dy.
\end{eqnarray}

\begin{lemma} \label{lema frobenius}
	If  the system (\ref{dA}-\ref{dY}) is compatible (i.e. has solutions) then 
$$
\begin{array}{l}
F_x=3FY+\dfrac{L_x}{L}F+LK_x, \\
\\
	F_y=3FR+\dfrac{L_y}{L}F-LK_y,  \\
	\\
	F^2=\dfrac{1}{3}LK_{xy}-\dfrac{1}{2}LK_x R-\dfrac{1}{2}LK_yY.
	\end{array}
$$
\end{lemma}
{\it Proof:} 
The condition $d(dY)=0$ gives 
$$
	F_x-3FY+\left(\dfrac{LL_{xy}-L_xL_y}{L^2}\right)Y+LKY-\dfrac{L_x}{L}F-LK_x=0,
$$
and after using (\ref{curvatura gaussiana})
\begin{equation}\label{Fx}
	F_x=3FY+\dfrac{L_x}{L}F+LK_x.
\end{equation}
Similarly from  $d(dR)=0$ one gets 
\begin{equation}\label{Fy}
	F_y=3FR+\dfrac{L_y}{L}F-LK_y.
\end{equation}
Thus the forms for  $dR$ and $dY$ are closed if and only if 
\begin{equation} \label{dF}
dF=\left(3FY+\dfrac{L_x}{L}F+LK_x\right)dx+\left(3FR+\dfrac{L_y}{L}F-LK_y\right)dy.
\end{equation}
Finally  $d(dF)=0$ results in 
\begin{equation}\label{F}
	F^2=\dfrac{1}{3}LK_{xy}-\dfrac{1}{2}LK_x R-\dfrac{1}{2}LK_yY. 		
\end{equation}
\hfill $\Box$\\
Now if we define $F$  in (\ref{dR},\ref{dY}) by (\ref{F}) we get a system of finite type: all the derivatives of $R,Y$ are given.  
\begin{corollary}\label{corinv}
The system (\ref{dR},\ref{dY}) with $F$ given by (\ref{F}) is in involution if and only if holds (\ref{Fx}) and (\ref{Fy}) due to (\ref{dR},\ref{dY},\ref{F})   
\end{corollary}
If the system (\ref{dR},\ref{dY}) with $F$ given by (\ref{F}) is in involution then  system (\ref{dA}-\ref{dY}) with $F$ given by (\ref{F}) is also in  involution and,  by the Frobenius Theorem, the initial values   $A_0, B_0, C_0, R_0, Y_0$ at point  $(x_0,y_0)$ define uniquely a (germ of) solution to the system (\ref{dA}-\ref{dY}). Thus one can think of  $A_0, B_0, C_0, R_0, Y_0$ as local coordinates on the space of integrals, even though this space may be not well defined globally.   

Note that the subsystem (\ref{dR},\ref{dY}) does not involve  $A, B$ and $C$. Therefore any solution to this subsystem fixes an orbit of $PSL_2(\mathbb{R})$-action, since  at a point $(x_0,y_0)$ any two triples of initial values for  $A, B, C$ are   M\"obius equivalent. 

Given a solution $R$ e $Y$ to (\ref{dR},\ref{dY}), the system (\ref{dA} - \ref{dC}) is {\it automorphic}: any two solutions $A,B,C$ are Möbius equivalent and the system  (\ref{dR},\ref{dY}) is {\it resolving} for (\ref{dA}-\ref{dC}), i.e. it describes the Möbius orbits on the space of solutions to (\ref{dA}-\ref{dC}) (see  \cite{O-82,V-04} for more detail).

\begin{lemma}\label{obs1}
	Suppose that the geodesic flow on a surface with metric (\ref{complex}) admits a fractional linear integral (\ref{fl}). Then the gaussian curvature is constant if and only if $$F=\dfrac{R_x-Y_y}{2}=0.$$
\end{lemma}
{Proof:} If the gaussian curvature $K$ is constant then 
$K_{xy}=K_x=K_y=0$ and therefore $F=0$ by  (\ref{F}). 
	
On the other hand, if $F=0$, then also $F_x=F_y=0$ and  equations (\ref{Fx}) and (\ref{Fy}) imply  
$LK_x=0,$ $LK_y=0$, hence
	 $K_x=K_y=0$ and  $K$ is constant.
\hfill $\Box$\\
\begin{lemma}\label{involut}
If $K=const$ then the system  (\ref{dR},\ref{dY}) is in involution. 
\end{lemma}
{\it Proof:} For $K=const$ we have $F=0$ and the claim follows by Corollary \ref{corinv} .  
\hfill $\Box$\\

For further analysis we need the following claim. 
\begin{proposition} \label{prop1}
	Suppose that the geodesic flow on a surface admits a fractional linear integral. Then if one partial derivative of the gaussian curvature vanishes then  the gaussian curvature is constant. 
\end{proposition}
{\it Proof:} 
	We can suppose that  $K_x=0 $ and  $K_y \neq 0$. Then from  equation  (\ref{F}) we have $Y=-\dfrac{2F^2}{LK_y}$. Differentiating by $x$ we get
	$$
		Y_x = -\dfrac{12F^2}{LK_y}Y-\dfrac{2L_xF^2}{L^2K_y}.
$$
Comparing with  (\ref{dY}) we obtain
	$$L^2K_yY^2 + (12LF^2 + LL_xK_y)Y+ 2F^2L_x=0.$$ Substituting $Y=-\dfrac{2F^2}{LK_y}$ results in 
$
		\dfrac{20F^4}{K_y}=0.
$
Hence $F=0$ and by Lemma  \ref{obs1} the curvature $K$ is constant.
\hfill $\Box$\\

\begin{theorem} \label{teo K não constante}
	Suppose that the geodesic flow on a surface  admits a fractional linear integral (\ref{fl}). Then the dimension of the space of such integrals is 3 if and only if the gaussian curvature is not constant.
\end{theorem}
{\it Proof:}
	Suppose that the gaussian curvature $K$ is not constant. Then $F \neq 0$ and by  Proposition \ref{prop1} holds $K_x \neq 0$ and $K_y \neq 0$. 
Let us solve equation    (\ref{F}) for $R$
\begin{equation}\label{R}
R=\dfrac{2K_{xy}}{3K_x}-\dfrac{K_y}{K_x}Y-\dfrac{2F^2}{LK_x}
\end{equation}
and consider the system of 2 equations (\ref{dY},\ref{dF}) with $R$ given by (\ref{R}). Then this system for $F,Y$ is of finite type  and its compatibility conditions reads as (\ref{dR}).

	Lets us show that $F$ and  $Y$ are fixed.

Differentiating (\ref{R}) we determine 
$$	
\begin{array}{l}
R_x=\left(\dfrac{K_{xx}K_y}{K_x^2} -\dfrac{K_{xy}}{K_x}\right)Y-\dfrac{K_{y}}{K_x}Y_x+\left(\dfrac{2L_x}{L^2K_x}+\dfrac{2K_{xx}}{LK_x^2}\right)F^2 -\dfrac{4F_xF}{LK_x}+\dfrac{2K_{xxy}}{3K_x}-\dfrac{2K_{xx}K_{xy}}{3K_x^2},\\  \\
R_y=\left(\dfrac{K_{xy}K_y}{K_x^2} -\dfrac{K_{yy}}{K_x}+\right)Y-\dfrac{K_{y}}{K_x}Y_y+\left(\dfrac{2L_y}{L^2K_x}+\dfrac{2K_{xy}}{LK_x^2}\right)F^2 -\dfrac{4F_yF}{LK_x}+\dfrac{2K_{xyy}}{3K_x}-\dfrac{2K_{xy}^2}{3K_x^2}. 
	\end{array}
	$$
Comparing with (\ref{dR}) and substituting the above expression for $R$ , we get 
	\begin{eqnarray} \label{eq26}
	\left(\dfrac{2L_x}{L^2K_x}-\dfrac{2K_{xx}}{LK_x^2}\right)F^2+\left(\dfrac{10F^2}{LK_x}+a_1\right)Y+5F-a_2=0,
	\end{eqnarray}
	where
	\begin{eqnarray*}
		&&a_1=-\dfrac{K_{xx}K_y}{K_x^2}+\dfrac{5K_{xy}}{3K_x}+\dfrac{L_xK_y}{LK_x};\\ \nonumber \\
		&&a_2=\dfrac{2K_{xxy}}{3K_x}-\dfrac{2K_{xx}K_{xy}}{3K_x^2}+\dfrac{L_{xy}}{L}-\dfrac{L_xL_y}{L^2}.
	\end{eqnarray*}
and 
	\begin{eqnarray} \label{eq27}
		\frac{10K_{xy}}{3LK_x^2}F^2-\dfrac{20}{L^2K_x^2}F^4 -\left(\dfrac{10K_yF^2}{LK_x^2}+b_1\right)Y-\dfrac{5K_y}{K_x}F+ b_2=0,
	\end{eqnarray}
	where
	\begin{eqnarray*}
		&&b_1=\dfrac{5K_{xy}K_y}{3K_x^2}-\dfrac{K_{yy}}{K_x}+\dfrac{L_yK_y}{LK_x},\\ \nonumber \\
		&&b_2=-\dfrac{2K_{xyy}}{3K_x}+\dfrac{10K_{xy}^2}{9K_x^2}+\dfrac{2L_yK_{xy}}{3LK_x}-\dfrac{L_{xy}K_y}{LK_x}+\dfrac{L_xL_yK_y}{L^2K_x}.
	\end{eqnarray*}
	
If $\dfrac{10F^2}{LK_x}+a_1 \neq 0$ then we can resolve  equation  (\ref{eq26}) for  $Y$. As the equation  (\ref{eq27}) has a non-vanishing coefficient of $F^4$, substituting the found expression for  $Y$ into (\ref{eq27}) we get a nontrivial equation for $F$. Thus $F$ is fixed and consequently $Y$  is  eventually fixed  by (\ref{eq26}).
		 
If $\dfrac{10F^2}{LK_x}+a_1 = 0$ then $F^2=-\dfrac{a_1LK_x}{10}$ and  $F$ is fixed as well as $F_x$. Finally, from  equation  (\ref{Fx}) we see that  $Y$ is fixed.  Thus, locally  there is only one  orbit of the Möbius group and the dimension is 3.
		
	The converse claim follows from Lemma \ref{involut}.
\hfill $\Box$
\begin{theorem} \label{teo1}
	Suppose that the geodesic flow on a surface admits a fractional linear integral (\ref{fl}). Then the dimension of the space of linear fractional integrals is 5 if and only if the gaussian curvature is constant.
\end{theorem} \label{teo K constante}
{\it Proof:} 
Follows from Lemma \ref{involut} and Theorem \ref{teo K não constante}. 
\hfill $\Box$

\section{Geodesic Nakai's 4-webs}\label{Geometria}

In this section we consider geometric questions  and therefore return to conformal coordinates for the metric   $g=\Lambda(x,y)\left(dx^2+dy^2\right)$. 

Suppose there is a fractional-linear integral (\ref{fl}) of the geodesic flow. 
It can be written as a function of local coordinates and the tangent vector 
 $v=(\xi,\eta)$. By the canonical isomorphism $\psi:TM \rightarrow T^{*}M$, the momentum $(p,q)$ and the tangent vector $(\xi,\eta)$ are related as follows: $$(\xi,\eta)=\left(\dfrac{p}{\Lambda},\dfrac{q}{\Lambda}\right).$$ 
So, in conformal coordinates, the integral is 
$$
	I=\dfrac{A\xi+B\eta}{C\xi+E\eta}.
$$
Given such integral $I$, we construct a one-parameter family of foliations as follows: for any  $\lambda \in \mathbb{R}$, the equation 
$$
		I=\dfrac{A\xi+B\eta}{C\xi+E\eta}=\lambda,
$$
gives an ODE
$$
(A-C\lambda)+(B-E\lambda)\dfrac{dy}{dx}=0
$$ on our surface after 
 $\xi=dx$ e $\eta=dy$.  
Integral curves of this ODE form a foliation on the surface. 
\begin{proposition}\label{geofol}
	The constructed foliation is geodesic.
\end{proposition}
{\it Proof:}
Follows from the uniqueness of an integral curve passing through $(x_0,y_0)$: the function $I$ is constant and is equal to $\lambda$ along the geodesic with the direction $(\xi,\eta)$ at $(x_0,y_0)$ fixed by $I=\lambda$. 
\hfill $\Box$\\

\begin{definition}
Nakai's 4-web is a planar 4-web such that the cross-ratio of the four directions tangent to the web leaves is constant.
\end{definition} 
\begin{theorem}\label{Teo Geo}
There exists a fractional linear integral (\ref{fl}) of the geodesic flow on a surface if and only if the surface carries a geodesic Nakai's 4-web. 
\end{theorem}
{\it Proof:}  Suppose there is an integral. The equation defining the geodesic foliations of Proposition \ref{geofol} can be written as 
	\begin{equation}\label{eq da ação}
		\dfrac{A+BP}{C+EP}=\lambda,
	\end{equation}
	where  $P=\dfrac{dy}{dx}$ is the inclination of the foliation  leaf. The integral gives the map $m: S\to SL_2(\mathbb{R})$,  $$m(x,y)=\left(\begin{array}{l r}
		A & B\\
		C & E
	\end{array}\right).$$ Then  $m\cdot P=\lambda$ and  $P=m^{-1}\cdot\lambda$.
	
	Consider a $4$-web, whose foliations are fixed by four values  $\lambda_i$, $i=1,2,3,4$. Then the cross-ratio of the inclinations $P_i$ is equal to the cross-ratio of $\lambda$-s and therefore is constant.

Now suppose that there are 4 geodesic  foliations tangent to 4 direction fields $\partial_x+J\partial_y$, $\partial_x+M\partial_y$, $\partial_x+N\partial_y$ e $\partial_x+T\partial_y$, such that
	\begin{equation}\label{razao cruz}
		\dfrac{J-M}{M-N}\cdot\dfrac{N-T}{T-J}=r
	\end{equation} with constant $r$. 
Resolving for $T$ we get
	\begin{equation*}
		T=\dfrac{rJM-(r-1)JN-MN}{J+(r-1)M-rN}
	\end{equation*}
and 
	\begin{eqnarray*}
	T_{x}=\dfrac{r(r-1)(M-N)^2J_x+r(J-N)^2M_x-(r-1)(J-M)^2N_x}{(J+(r-1)M-rN)^2},
	\end{eqnarray*}
	\begin{eqnarray*}
			T_{y}=\dfrac{r(r-1)(M-N)^2J_y+r(J-N)^2M_y-(r-1)(J-M)^2N_y}{(J+(r-1)M-rN)^2}.
	\end{eqnarray*}
The  total differentiation by $x$ along an integral curve  $y=y(x)$ of the field $\partial_x+P\partial_y$ is 
	\begin{equation}\label{deriv total}
		\dfrac{d^2y}{dx^2}=P_x+PP_y.
	\end{equation}
If these curves are geodesics then $\partial_x+P\partial_y$ is  the Jacobi field.  Using the equation	
$$
		\dfrac{d^2y}{dx^2}=-\Gamma_{11}^2+(\Gamma_{11}^1-2\Gamma_{12}^2)\dfrac{dy}{dx}-(\Gamma_{22}^2-2\Gamma_{12}^1)\left(\dfrac{dy}{dx}\right)^2+\Gamma_{22}^1\left(\dfrac{dy}{dx}\right)^3,
$$
for unparametrized geodesics, we get 
	\begin{equation*}
		P_x+PP_y=-\Gamma_{11}^2+(\Gamma_{11}^1-2\Gamma_{12}^2)P-(\Gamma_{22}^2-2\Gamma_{12}^1)P^2+\Gamma_{22}^1P^3,
	\end{equation*}
where $\Gamma_{jk}^i$ are Christoffel symbols for  the Levi-Civita connection.
Therefore for the fields $\partial_x+J\partial_y$, $\partial_x+M\partial_y$, $\partial_x+N\partial_y$ e $\partial_x+T\partial_y$ we have
	\begin{eqnarray}\label{eq geo1}
	&&	J_{x}+JJ_{y}-\frac{\Lambda_y}{2\Lambda}+\frac{\Lambda_x}{2\Lambda}J-\frac{\Lambda_y}{2\Lambda}{J}^2+\frac{\Lambda_x}{2\Lambda}{J}^3=0,\\ \label{eq geo2}\nonumber \\
	&&	M_{x}+MM_{y}-\frac{\Lambda_y}{2\Lambda}+\frac{\Lambda_x}{2\Lambda}M-\frac{\Lambda_y}{2\Lambda}{M}^2+\frac{\Lambda_x}{2\Lambda}{M}^3=0,\\ \label{eq geo3}\nonumber \\
	&&	N_{x}+NN_{y}-\frac{\Lambda_y}{2\Lambda}+\frac{\Lambda_x}{2\Lambda}N-\frac{\Lambda_y}{2\Lambda}{N}^2+\frac{\Lambda_x}{2\Lambda}{N}^3=0,\\ \label{eq geo4}\nonumber \\
	&&	T_{x}+TT_{y}-\frac{\Lambda_y}{2\Lambda}+\frac{\Lambda_x}{2\Lambda}T-\frac{\Lambda_y}{2\Lambda}{T}^2+\frac{\Lambda_x}{2\Lambda}{T}^3=0.
\end{eqnarray}
Resolving (\ref{eq geo1},\ref{eq geo2},\ref{eq geo3})  for $J_x, M_x, N_x$,  we get 
\begin{eqnarray}
	&&J_{x} = -\dfrac{{\Lambda_xJ}^3- \Lambda_y{J}^2 + (2\Lambda{J}_y + \Lambda_x)J - \Lambda_y}{2\Lambda}, \label{Jx}\\ \nonumber \\
		&&M_{x}= -\dfrac{\Lambda_x{M}^3- \Lambda_y{M}^2 + (2\Lambda{M}_y + \Lambda_x)M- \Lambda_y}{2\Lambda}, \label{Mx}\\ \nonumber \\
		&&N_{x}= -\dfrac{\Lambda_x{N}^3-\Lambda_y{N}^2 +(2\Lambda{N}_y + \Lambda_x)N- \Lambda_y}{2\Lambda} \label{Nx}.
\end{eqnarray}
	
Substituting $J_{x}, M_{x}, N_{x}, T, T_{x}$, and  $T_{y}$ into  (\ref{eq geo4}), we obtain
	\begin{eqnarray}\label{eq 4 boa}
	\dfrac{r(r - 1)(J - M)(J - N)(M - N) u}{2\Lambda(J+(r-1)M-rN)^3}=0,
	\end{eqnarray}\normalsize
where 
$$
	u=\Lambda_x(J-M)(J-N)(M-N)+2\Lambda((J-M)N_y-(J-N)M_y+(M-N)J_y).
$$
Each factor in the product $r(r - 1)(J - M)(J - N)(M - N)$ is not zero hence $u=0$. Note also that $J+(r-1)M-rN\neq 0$, since otherwise  from (\ref{razao cruz}) follows $J=N$.
	
Now we construct a map  
	$m: S \rightarrow SL_2(\mathbb{R})$ , 
$m=\left(\begin{array}{cc} A & B  \\ C & E\end{array}\right),$ resolving 
\begin{eqnarray} \label{eq A B C E}
\dfrac{A+BJ}{C+EJ} = 0, \ \ \
\dfrac{A+BM}{C+EM} = 1, \ \ \
C+EN=0, \ \ \
AE-BC=1, 
\end{eqnarray}
for $A,B,C,E$, and show that 	$I=\dfrac{Ap+Bq}{Cp+Eq}$ is a first integral. 
To this end we verify that 
\begin{equation}\label{eqintnorm} 
\begin{array}{l}
	2\Lambda( CA_x-AC_x) - \Lambda_y=0,  \\
	\\ 
	2\Lambda(EA_x-AE_x+CA_y-AC_y+CB_x-BC_x) + \Lambda_x=0,\\
	\\
	2 \Lambda (EA_y-AE_y+CB_y-BC_y+EB_x-BE_x) - \Lambda_y=0, \\
	\\
	2\Lambda(EB_y-BE_y) + \Lambda_x=0. 
\end{array}
\end{equation}
Then $I$ is a first integral by Lemma  \ref{lema sistema homo}. 
	
From (\ref{eq A B C E}) we obtain
$$
A = \dfrac{J(M-N)E}{J-M},\ \ \ 
B= -\dfrac{(M-N)E}{J-M},  \ \ \ 
C= - NE.
$$
Differentiating, we have 
$$
	\begin{array}{l}
	A_x=-\dfrac{EM(M-N)J_x}{(J-M)^2} + \dfrac{EJ(J-N)M_x}{(J-M)^2} + \dfrac{E_xJ(M - N)}{J - M}- \dfrac{EJN_x}{J - M},\\   \\
	A_y=-\dfrac{EM(M-N)J_y}{(J-M)^2} + \dfrac{EJ(J-N)M_y}{(J-M)^2} + \dfrac{E_yJ(M - N)}{J - M}- \dfrac{EJN_y}{J - M}, 
	\end{array}
$$
as well as 
$$
\begin{array}{l}
B_x=\dfrac{E(M - N)J_x}{(J - M)^2} - \dfrac{E(J - N)M_x}{(J - M)^2} - \dfrac{E_x(M - N)}{J - M}+ \dfrac{EN_x}{J - M},\\ 
 \\
B_y=\dfrac{E(M - N)J_y}{(J - M)^2} - \dfrac{E(J - N)M_y}{(J - M)^2}  - \dfrac{E_y(M - N)}{J - M}+ \dfrac{EN_y}{J - M},
\end{array}
$$
and
$$
\begin{array}{l}
C_x= -EN_x - E_xN, \\ 
C_y= -EN_y - E_yN. 
\end{array}
$$
Now the expression  $2\Lambda( CA_x-AC_x) - \Lambda_y(AE-BC)$ writes as 
\begin{eqnarray}
\dfrac{2\Lambda E^2MN(M - N)J_x}{(J - M)^2}- \dfrac{2\Lambda E^2JN(J - N)M_x}{(J - M)^2} + \dfrac{2\Lambda E^2JMN_x}{J - M} - \dfrac{\Lambda_yE^2(M - N)(J - N)}{J - M}.\nonumber
\end{eqnarray}
One checks that 
$$
	E^2=\dfrac{J-M}{(M-N)(J-N)}, 
$$
and has 	
$$
	2\Lambda( CA_x-AC_x) - \Lambda_y(AE-BC)=-\dfrac{uJMN}{(J-M)(M-N)(J-N)}.
$$
Similarly we get
$$
 	2\Lambda(EA_x-AE_x+CA_y-AC_y+CB_x-BC_x) + \Lambda_x(AE-BC)= \dfrac{u(JN +JM +MN)}{(J-M)(M-N)(J-N)}
$$
and 
$$
	2\Lambda (EA_y-AE_y+CB_y-BC_y+EB_x-BE_x) - \Lambda_y(AE-BC)=-\dfrac{u(J+M+N)}{(J-M)(M-N)(J-N)},
$$
and 
$$
2\Lambda(EB_y-BE_y) + \Lambda_x(AE-BC)=\dfrac{u}{(J-M)(M-N)(J-N)}.
$$
Since $u=0$ and $AE-BC=1$, the last 4 equations give (\ref{eqintnorm}).
\hfill $\Box$\\

\section{Concluding remarks}
One easily checks that the explicit examples of metrics constructed by Agapov and Shubin \cite{AS-21} admit {\it projective} vector fields, whose local flows maps geodesics into geodesics but do not have to respect metrics (see \cite{L-82}). It would be interesting to understand relation between existence conditions for projective vector fields and for fractional-linear integrals. 

\section*{Acknowledgments}
This research was  supported by FAPESP grant \# 2022/12813-5 (S.I.A) and CAPES grant \# 88882.434346/2019-01 (T.G.P.A).

\end{document}